\documentclass[reqno]{amsart}
\newtheorem{theorem}{ Theorem}[section]
\newtheorem{remark}{ Remark}[section]
\def\sign{\mathop{\rm sgn}\nolimits}
\begin{document}
\title [The non-viscous Burgers equation associated with random positions]
{The non-viscous Burgers equation associated with random positions
in coordinate space: a threshold  for blow up behaviour}

\author[Albeverio,\,Rozanova]{Sergio Albeverio $^{1}$, Olga Rozanova $^{2}$}

\address[$^{1}$]{Universit\"{a}t Bonn,
Institut f\"{u}r Angewandte Mathematik, Abteilung f\"{u}r
Stochastik, Wegelerstra\ss e 6, D-53115, Bonn; IZKS, Bonn; BiBoS,
Bielefeld--Bonn}
\address[$^{2}$]{Mathematics and Mechanics Faculty, Moscow State University, Moscow
119992, Russia}

\thanks {Supported by  DFG 436 RUS 113/823/0-1.}

\email[$^{1}$]{albeverio@uni-bonn.de}
\email[$^{2}$]{rozanova@mech.math.msu.su}


%
%


\date{\today}


\begin{abstract}
It is well known that the solutions to the non-viscous Burgers
equation develop a gradient catastrophe at a critical time provided
the initial data have a negative derivative in certain points. We
consider this equation assuming that the particle paths in the
medium are governed by a random process with a variance which
depends in a polynomial way on the velocity. Given an initial
distribution of the particles which is uniform in space and with the
initial velocity linearly depending on the position we show both
analytically and numerically that there exists a threshold effect:
if the power in the above variance is less than 1, then the noise
does not influence the solution behavior, in the following sense:
the mean of the velocity when we keep the value of position fixed
goes to infinity outside the origin. If however the power is larger
or equal 1, then this mean decays to zero as the time tends to a
critical value.
\end{abstract}
\maketitle

\medskip


\section{Introduction}

The non-viscous Burgers equation is perhaps the simplest equation
that models the nonlinear phenomena in a force free mass transfer.
It has the form
$$
u_t+(u,\nabla)\,u=-\beta u,\eqno(1.1)
$$
where $u(x,t)=(u_1,...,u_n)$ is a vector-function ${\mathbb
R}^{n+1}\to{\mathbb R}^n,$ $\beta\ge 0$ is a constant friction
coefficient. Consider the Cauchy data
$$u(x,0)=u_0(x).\eqno(1.2)$$
Problem (1.1), (1.2) has an implicit solution $$u(t,x)=e^{-\beta
t}\,u_0(x-\frac{1}{\beta}\,(e^{\beta t}-1) u(t,x)),$$ for $\beta>0$
and $$u(t,x)=u_0(x-t u(t,x)),$$ for $\beta=0.$

In several cases we can obtain an explicit solution. For example, if
$$u_0(x)=\alpha x,\quad \alpha\in{\mathbb R}, \eqno(1.3)$$ one easily
gets
$$
u(t,x)=\frac{\alpha x e^{-\beta
t}}{1+\frac{\alpha}{\beta}(1-e^{-\beta t})},\quad \beta>0,
\eqno(1.4)
$$
and, for $\, \beta=0:$
$$u(t,x)=\frac{\alpha
x}{1+\alpha t} . \eqno(1.5)
$$ Thus, if $\alpha<-\beta,$
the solution develops a singularity at the origin as $t\to
T,\,0<T<\infty,$ where
$$T=\frac{1}{\beta}\,\ln
\frac{\alpha}{\alpha+\beta},\quad \mbox{for }\,\beta>0,\qquad
T=-\frac{1}{\alpha},\quad \mbox{for }\,\beta=0,\quad \alpha\ne
0.\eqno(1.6)$$ This phenomenon is called the gradient catastrophe.
It is well known (see \cite{Hopf})that a viscous perturbation of
form $\sigma \Delta u,\,\sigma>0,$ entails a globally in time smooth
solution (at least for bounded initial data). An exceptional case is
exactly given by a solution which is linear in $x$ as mentioned
above, which does not feel the viscous term.

Our main question is: can a stochastic perturbation  suppress the
appearance of unbounded gradients?

We can introduce the Lagrangian coordinate $x(t)$ to label a point
which moves together with the medium, that is $\frac{d
x(t)}{dt}=u(t,x(t)):=u_1(t).$ Thus, $x=x(t)$ is the equation for the
particle path, when the particle moves along the Burgers fluid.
Equation (1.1) is equivalent to the following system of ODE:
$$
\dot x(t)= u_1(t),\quad \dot  u_1(t)=-\beta  u_1 (t)\eqno(1.7)
$$
Further on we will omit the index 1.

In the theory of stochastic dynamical systems one often considers a
stochastic perturbation of the velocity, which leads to the
appearance of a white noise in the second of equations (1.7). The
problem of solving such equations was investigated in many works
(see, e.g. \cite{AKl1}, \cite{AK2}, \cite{AK}, \cite{AHZ},
\cite{Risken}). This type of stochastic perturbation corresponds to
the stochastically forced Burgers equation, or in the language of
physicists, Burgers turbulence. This has been an area of intensive
research activity in the last decade (see e.g. \cite{Woy}, and for a
very recent review \cite{Khanin}, and references therein).

 The behavior of the gradient of velocity was studied
earlier in other contexts in \cite{Bouchaud}, \cite{Gurarie}, but
this problem is quite different from the problem considered in this
paper.

 In the present paper we consider a medium with random particles
paths, more precisely, described by a $2\times n$ dimensional
It$\rm\hat{o}$ stochastic differential system of equations
$$
d X_k(t)=U_k(t)\,dt +\sigma |U(t)|^p\,d(W_k)_t, $$
$$d U_k(t)=\,-\beta\,U_k(t)\,dt,\quad k=1,..,n, \eqno(1.8)
$$
$$
X(0)=x,\quad U(0)=u,\quad t\ge 0,
$$
where $(X(t),U(t))$ runs in the phase space ${\mathbb R}^n\times
{\mathbb R}^n,$ \,$\sigma>0$ and $p\ge 0$ are constants,
$(W)_t=(W)_{k,t},\,k=1,...,n,$ is the $n$ - dimensional Brownian
motion. We remark that for $p\le 1$ for initial distributions of $x$
and $u$ from the class $L^2({\mathbb R}^n)$ one can guarantee a
global existence of a unique solution to (1.8)\cite{oksendal}.

Let us denote by $\hat u(t,x)$ the mean of the velocity $U(t)$ at
time $t$ when we keep the value of $X(t)$ at time $t$ fixed but
allow $U(t)$ to take any value it wants (e.g.\cite{Chorin}) and
chose this function for comparison with the solution to the solution
of the non-viscous Burgers equation.

We can interpret system (1.8) also as follows: assume that we
measure the position of a particle with an error depending on its
velocity and then try to restore the velocity. If the coefficient
 $p$ increases, the error for large velocities increases, too. It is natural
 to expect that the mean of difference
 between two "very indefinite" neighbor coordinates necessary to calculate the velocity
 tends to zero.

Can we hope to extract from our measurement a realistic information
on  such critical phenomena as the blow up occurring in a medium
described by the Burgers equation associated with (1.8)? As we will
see, at least for uniform initial distribution of particles (in the
sense of Sec.3, see below) the answer depends on the exponent $p,$
namely, if $p\ge 1,$ the information gets lost. The threshold value
$p=1$ is not astonishing. The sub-linear rase of the drift and
diffusion ($p\le 1$) coefficients warrants the global existence of
the SPDE solution provided the initial distributions of the particle
positions and the velocities are square integrable. Therefore one
can hope that the solutions to the SPDE with a sub-linear diffusion
coefficient behave  in a some "predictable" way.

For example, for $\beta=0$ the function $\hat u(t,x)$ demonstrates
the following behaviour near the critical time $T$. For $p\in [0,1)$
$$\hat u(t,x)\, = \,
\frac{\alpha}{1-\frac{t}{T}}\,x\,+o\left(\frac{1}{1-\frac{t}{T}}\right),\quad
t\to T,\quad x\in{\mathbb R}^n;$$ for $p>1$
$$
\hat u (t,x)\,=\,\,-\,C\,|x|^{\frac{2(1-p)}{p}}\,x
\left(1-\frac{t}{T}\right)\,+\,o(1-\frac{t}{T}),\quad t\to T,
\,\,x\ne 0,$$ where $C$ is the positive constant depending in
particular on the dimension of space, however,
$$\hat u (t,x)\,=\, \frac{\alpha}{1-\frac{t}{T}}\,x\,+\,o\left(|x|\right),\quad  t\in
[0,T),\,\, x\to 0.$$ In other words, for $p\in [0,1)$ we see that
$\hat u(t,x)\sim u(t,x),\, t\to T$ at every point $ x\in{\mathbb
R}^n,$ where $u(t,x)$ is given by (1.5), i.e. $\hat u(t,x)$ keeps
the property of solutions to the inviscid Burgers to have a gradient
catastrophe. For $p\ge 1$ outside the origin $x=0$ the function
$\hat u (t,x)\to 0,\,t\to T,$ and the features of the solution of
the non-perturbed equation fail.  A jump is being formed at the
origin $x=0$ but the height of this jump at the same time tends to
zero as $t\to T.$

The paper is organized as follows. In Sec.2 we find an exact
solution to the Fokker-Planck equation with special initial data and
derive an integral formula for $\hat u.$ In Sec.3 we consider the
simplest case $p=0,$ where for several initial distribution it is
possible to find $\hat u$ exactly. In Sec.4 we formulate the main
theorem for the case of uniform initial distribution concerning
asymptotics of $\hat u$ near the critical time $T$ and the origin
$x=0$ and prove it. Here we  present also the results of
computations given directly according to the integral formula for
$\hat u$ that illustrate the fidelity of our asymptotic formulas. In
Sec.5 we present certain results, both asymptotic and numeric,
concerning a Gaussian initial distribution and show that for $0\le p
< 1$ the behaviour of $\hat u$ is close to the behaviour of the same
function at $p=0,$ where the exact formula can be derived. In Sec.6
we discuss the question on the "observable" density and "induced"
velocity. In Conclusion we sum up our results and discuss the
question on a possibility to construct a PDE such that $\hat u$ is
its solution. We also argue possible applications of the results.

\section{Exact solution to the Fokker-Planck equation for special initial data}

The Fokker-Planck equation associated to (1.8) for the probability
density in position and velocity space $P= P(t,x,u)$ has the form
$$
\frac{\partial P(t,x,u)}{\partial t
}=\left[-\sum\limits_{k=1}^n\,u_k\,\frac{\partial }{\partial x_k }\,
+\beta \,\sum\limits_{k=1}^n\, \left(u_k \,\frac{\partial }{\partial
u_k }\,+\, 1 \right)\, + \frac{1}{2}\, \sigma^2
|u|^{2p}\,\frac{\partial^2 }{\partial x_k^2}\,\right] P(t,x,u),
\eqno(2.1)
$$
subject to the initial data
$$P(0,x,u)=P_0(x,u).$$

Let us set $\Omega_L:=[-L,L]^n,\,L>0.$ Thus,
$$
\hat
u(t,x)=\,\lim\limits_{L\to\infty}\,\frac{\int\limits_{\Omega_L}\,u\,P(t,x,u)\,du}{\int\limits_{\Omega_L}\,P(t,x,u)\,du},\quad
t\ge 0,\,x\in \Omega_L,\eqno(2.2)
$$
provided the limit exists.

If we choose
$$
P_0(x,u)=\delta (u-u_0(x))\,f(x)=\prod\limits_{k=1}^{n}\,\delta
(u_k-(u_0(x))_k)\,f(x),\eqno(2.3)$$ with an arbitrary sufficiently
regular $f(x),$ then
$$
\hat u(0,x)=u_0(x).
$$
The function $f(x)$ has the meaning  of a probability density of the
particle positions in the space at the initial moment of time.

Moreover, formula (2.2) can be investigated for functions $f(x)$
which are not necessarily probabilities densities.

Let us choose $$u_0(x)=\alpha x,\quad\alpha<0\eqno(2.4)$$ as initial
data of the non-perturbed Burgers equation. One can see from (1.4),
(1.5) that the gradient of the solution become unbounded as $t\to
T.$ Thus, we are interested in the behavior of $\hat u(t,x)$
comparing with the solution $u(t,x)$ to (1.1).

We apply formally the Fourier transform in the variable $x$ to
(2.1), (2.3) (2.4) to obtain for $\,\tilde P\,= \,\tilde
P(t,\lambda, u)\,$
$$
\frac{\partial  \tilde P}{\partial
t}\,=\,\beta\,\sum\limits_{k=1}^n\,u_k\, \,\frac{\partial \tilde
P}{\partial u_k}\,+\,(\beta - \frac{\sigma^2}{2} \,|u|^{2p}
\,|\lambda|^2\,-i\,(\lambda ,u))\,\tilde P,\eqno(2.5)$$
$$
\tilde P(0,\lambda, u)=\frac{1}{(|\alpha| \sqrt{2
\pi})^n}\,e^{-i\frac{(\lambda,u)}{\alpha}}\,f\left(\frac{u}{\alpha}\right).\eqno(2.6)
$$
Equation (2.5) is of the first order, therefore the Cauchy problem
(2.5), (2.6) for the function $\tilde P (t,\lambda, u)$ can easily
be solved. Thus, for $\beta=0$
$$
\tilde P (t,\lambda, u)=\frac{f\left({u}/{\alpha}\right)}{(|\alpha|
\sqrt{2 \pi} )^n}\, e^{-\frac{\sigma^2}{2}\, |u|^{2p} |\lambda|^2 t
-i(\lambda,u)(\frac{1}{\alpha}+t)},\eqno(2.7)
$$
for $\beta>0$
$$
\tilde P (t,\lambda, u)=\frac{f\left({u\,e^{\beta
t}}/{\alpha}\right)}{(|\alpha| \sqrt{2 \pi} )^n}\, e^{\beta
t-\frac{\sigma^2}{2 p \beta}\, |u|^{2p} |\lambda|^2\,(e^{\beta p
t}-1) -i(\lambda,u)((\frac{1}{\alpha}+\frac{1}{\beta})e^{\beta
t}-\frac{1}{\beta})}.\eqno(2.8)
$$
Further, the inverse Fourier transform gives for $\beta=0$
$$
P(t,x,u)=\frac{f\left({u}/{\alpha}\right)}{(|\alpha|
\sigma\sqrt{2\pi
t})^n\,|u|^{pn}}\,e^{-\frac{|u(\frac{1}{\alpha}+t)-x|^2}{2\sigma^2 t
|u|^{2p}}}, \quad t>0, \eqno(2.9)
$$
and for $\beta>0$
$$
P(t,x,u)=\frac{f\left({u\,e^{\beta t}}/{\alpha}\right)\,e^{-\beta
t}}{(|\alpha|\sigma\sqrt{2\pi (e^{2\beta p t}-1)/\beta
t})^n\,|u|^{pn}}\,e^{-\frac{2\beta p\,|u\left(e^{\beta t}
(\frac{1}{\alpha}+\frac{1}{\beta})-\frac{1}{\beta}\right)-x|^2}{2\sigma^2
|u|^{2p}\,(e^{2\beta p t}-1)}}, \, t>0. \eqno(2.10)
$$
It is easy to see that  the limit for $\beta\to 0$ of (2.10) gives
(2.9).

Now  we substitute $P(t,x,u)$ in  (2.2) to get an integral
representation of the quantity expectation  $\hat u(t,x).$ It is
important to note that for every even function $f(x)$  the upper
integral in (2.2) vanishes at $x=0$ or $t=T,$  hence $\hat u(t,0)=
\hat u(T,0)=0.$

However,  representation (2.2), (2.9) (or (2.10)) does not allow to
see the behaviour of $\hat u.$ Therefore we firstly consider the
situations where the integrals in (2.2) can be explicitly computed
(Sec.3) or analyze the asymptotics in certain important points
(Sec.4).

\section{Case of constant noise variance}

Firstly we set $p=0.$ In this section we concentrate on a special
choice of the initial probability density $f(x),$ allowing to obtain
an explicit formula for $\hat u(t,x)$ if the initial data have the
form (2.4).

Let us consider uniform initial distribution of particles as
follow:
$$f(x)=\left\{\begin{array}{ll}
f_L=\frac{1}{(2L)^n}\,,&\quad x\in \Omega_L;\\ 0,&\quad otherwise.
\end{array}\right.\eqno(3.1)$$

Applying (2.2), (2.9), (2.10) we readily calculate
$$\hat u(t,x)=\left\{\begin{array}{ll}
\frac{\alpha x e^{-\beta t}}{1+\frac{\alpha}{\beta}(1-e^{-\beta
t})},&\quad t<T; \\ 0,&\quad t=T
\end{array}\right.,\eqno(3.2)$$
for $\beta>0$ and
$$\hat u(t,x)=\left\{\begin{array}{ll}
\frac{\alpha x}{1+\alpha t},&\quad t<T \\ 0,&\quad t=T
\end{array}\right.,\eqno(3.3)$$
for $\beta=0$. We also see that (3.3) results in the limit $\beta\to
0$ from (3.2).

It is enough to compare these expressions with (1.4), (1.5) to  see
that the white noise with a constant variance $\,\sigma^2\,$ does
not influence the mean of the velocity. Thus, at any point
$(T,x),\,x\ne 0,$ the function $\hat u (t,x)$ has a discontinuity
(as for the case without noise).

We can also compute the  variance
$$
\hat v(t,x)\,:=\,\lim\limits_{L\to
\infty}\,\frac{\int\limits_{\Omega_L}\,(u-\hat
u)^2\,P(t,x,u)\,du}{\int\limits_{\Omega_L}\,P(t,x,u)\,du},\quad t\ge
0,\quad x\in {\mathbb R}^n. \eqno(3.4)
$$
We have
$$
\hat v(t,x)=
\frac{\sigma
t}{(t+\frac{1}{\alpha})^2}
,\quad \beta =0,
$$
$$
\hat v(t,x) =\frac{\sigma t}{\left(e^{\beta
t}\left(\frac{1}{\alpha}+\frac{1}{\beta}\right)-\frac{1}{\beta}\right)},\quad
\beta
>0.
$$
Thus, as $t\to T,$ the denominator of $\hat v(t,x)$ becomes zero and
the possible values of $U (t)$ at any $x$ dissipate over the space.

A behavior of the mean of velocity which contrasts with the previous
one can be obtained, e. g., for
$f(x)=\left(\frac{k}{\sqrt{\pi}}\right)^n\,\exp (-k^2 x^2),\,k>0.$
It is easy to compute that in this case for $\beta=0$ one has
$$
\hat u(t,x)=\frac{(1+\alpha t)\,\alpha \, x}{\alpha^2
\,t^2\,+\,2(k^2\sigma^2+\alpha)\,t\,+1},\eqno(3.5)
$$
$$
\hat v(t,x)=\frac{\sigma^2\alpha^2\, t}{\alpha^2
\,t^2\,+\,2(k^2\sigma^2+\alpha)\,t\,+1}.
$$
The denominator does not vanish for any fixed $x,\,\alpha,\,\sigma,$
so $\hat u (t,x)\to 0,\quad \hat v(t,x)\to 0,\quad t\to \infty,$ and
$\hat u(t,x)$ is continuous at any point $(T,x).$ The minimal value
of ${\rm div} \hat u(t,x)$ is attained at the time
$t_*=\frac{1}{\alpha}\left( \sqrt{\frac{2}{-\alpha}}\,\sigma k
-1\right)<T.$ Moreover, $t_*>0$ only if $ \sigma
k<\sqrt{\frac{2}{-\alpha}}.$

Computations made for some special classes of $f(x)$ allow to
suggest that a similar behavior is provided by the function $\hat
u(t,x)$ if $f(x)=f(|x|)$ and $\int \limits_{{\mathbb
R}^n}\,|x|^2\,f(|x|)\,dx\,<\,\infty.$ For example, for the class of
$f(|x|)=\frac{const}{(1+k^2 x^2)^s},\,n=1$ we have the analytical
result
$$\hat u(t,x)\sim \frac{\alpha^2}{(2(s-1)-1)k^2}\,x\,(1+\alpha t),\quad \mbox{as}\,\,x\to 0,\quad s\in {\mathbb N},\,
s\ge 2.  $$ A numerical study suggests that for $n=1$ as $x\to 0$
$$\mbox{for}\, s>1 \qquad \hat u(t,x)\sim {\rm const}\,x\,(1+\alpha
t),$$
$$\mbox{for}\,s\in [\frac{1}{2},1]\qquad \hat u(t,x)\sim {\rm
const}\,x\,,$$$$\mbox{for}\,s<\frac{1}{2}\qquad
 \hat u(t,x)\sim \frac{\rm const}{1+\alpha t}\,x\,.$$
Here the constants do not depend on $t.$



\section{Asymptotic behavior for uniform initial distribution}

We consider again $f(x)$ given in (3.1). If $p>0,$ formula (2.2)
does not allow to compute $\hat u (t,x)$ explicitly. Thus, we need
to extract from this formula an information that allows us to answer
the main question of this article: does the stochastic perturbation
suppress the singularity formation?

We will show that there exists a critical value of the parameter,
$p=1,$ such that for $p\le 1$ the mean of the velocity behaves very
closely to the velocity in the case $p=0$ (as in the non
stochastically perturbed case for $t<T$). In contrast, for $p>1,$
$\hat u (t,x)$ vanishes as $t\to T$.

We deal with the case $\beta=0$ to avoid cumbersome formulas,  the
results for $\beta>0$ will be qualitatively the same, we present
them in Remark 4.1.

\begin{theorem}
The mean $\hat u (t,x) $ of the random variable $U(t),$ given the
position $X(t),$ where $X(t), U(t)$ solve the SDE (1.8), provided
$\hat u (0,x)=\alpha x, \,\alpha<0, $
is given by formula (2.2), where $P(t,x,u)$ is given by (2.9),
(2.10).

If $\beta=0$ and initially the particles are distributed uniformly
in the sense (3.1), the asymptotic behaviour of $\hat u(t,x)$ for
$t\to T$ can be analyzed explicitly.

Namely, for any $p\in [0,1)$ the mean $\hat u (t,x),$ being equal to
zero   at any point $x\in{\mathbb R}^n,\, t=T,$ is discontinuous at
every such point if $ \,x\ne 0.$ More precisely,
\begin{itemize}
\item for $p=0$ the mean $\hat u(t,x)$ coincides with the solution to
the problem (1.1), (2.4) for $t< T;$

\item for $p\in (0,1)$ the asymptotics $$\hat u(t,x)\, = \,
\frac{\alpha}{1-\frac{t}{T}}\,x\,+o\left(\frac{1}{1-\frac{t}{T}}\right),\quad
t\to T,\quad x\in{\mathbb R}^n,$$ takes place.

\item For $p\ge 1$ at any $x\in {\mathbb R}^n,\,x\ne 0,$ $|\hat u
(t,x)|\to 0$ as $t\to T,$ however, $ {\rm div}_x \hat u (t,x)\to
\infty,\, x=0,\,t\to T.$ More precisely, for any $x\ne 0,$
$$
\hat u (t,x)\,=\,\,-\,C\,|x|^{\frac{2(1-p)}{p}}\,x
\left(1-\frac{t}{T}\right)\,+\,o(1-\frac{t}{T}),\quad t\to T.$$
where $C$ is the positive constant given in (4.4), and for any $t\in
[0,T)$
$$\hat u (t,x)\,=\, \frac{\alpha}{1-\frac{t}{T}}\,x\,+\,o\left(|x|\right),\quad x\to 0.$$

\end{itemize}
\end{theorem}

\medskip

\begin{proof}
Let us set $\varepsilon:=1-\frac{1}{T}\, t, $ where the critical
time $T$ is introduced in (1.6) ($\varepsilon\in (0,1]$). We will
write below $t(\varepsilon)= (1-\varepsilon)T.$

\medskip
Firstly we compute asymptotics of (2.2) near the critical time
$t=T\,(\varepsilon=0)$.


Let us fix $x\ne 0.$  We have the following expansion as
$\varepsilon\to 0,\,$$p>1$
$$\int\limits_{\Omega_L}\,u\,P(t(\varepsilon),x,u)\,du \,=\,f_L \,\frac{\omega_n}{\alpha}\left(\sqrt{-\frac{1}{2\sigma^2\pi \alpha}}\right)^n
\,\int \limits_0^{|\alpha L|}\, e^{\frac{\alpha |x|^2}{2\sigma^2
|u|^{2p}}}\,\frac{1}{|u|^{(n+2)p-n-1}}\,d|u|\,\,x\varepsilon +
o(\varepsilon).\eqno(4.1)$$
$$\int\limits_{\Omega_L}\,P(t(\varepsilon),x,u)\,du \,=\,f_L\,\omega_n\,\left(\sqrt{-\frac{1}{2\sigma^2\pi\alpha}}\right)^n\,\int
\limits_0^{|\alpha L|}\, e^{\frac{\alpha |x|^2}{2\sigma^2
|u|^{2p}}}\,\frac{1}{|u|^{pn-n+1}}\,d|u| +
O(\varepsilon),\eqno(4.2)$$ where $\omega_n$ is the area of surface
of ($n-1$) - dimensional sphere, $n\ge 2;$  $\omega_1=2.$ Recall
that for fixed $x$ these integrals are functions of $\varepsilon$
only.


Both integrals in (4.1) and (4.2) converge as $L \to \infty$ and can
be expressed through Gamma-functions, so according to (2.2) we get
$$\hat u(t,x)\,=\,\frac{ \frac{1}{\alpha}
\left(\sqrt{-\frac{1}{2\sigma^2\pi\alpha}}\right)^n \,\int
\limits_0^\infty\, e^{\frac{\alpha |x|^2}{2\sigma^2
|u|^{2p}}}\,\frac{1}{|u|^{(n+2)p-n-1}}\,d|u|\,\,x\varepsilon +
o(\varepsilon)}{\left(\sqrt{-\frac{1}{2\sigma^2\pi\alpha}}\right)^n\,\int
\limits_0^\infty\, e^{\frac{\alpha |x|^2}{2\sigma^2
|u|^{2p}}}\,\frac{1}{|u|^{pn-n+1}}\,d|u| + O(\varepsilon)}=$$
$$\frac{\varepsilon\,x\,|x|^{(n+2)\frac{1-p}{p}}\,\frac{1}{\alpha}
\left(\sqrt{-\frac{1}{2\sigma^2\pi
\alpha}}\right)^n\,\left(\sqrt{-\frac{\alpha}{2\sigma}}\right)^{(n+2)\frac{1-p}{2p}}\,
\Gamma\Big((n+2)(p-1)/2p\Big)\,+o(\varepsilon)}
{|x|^{n(1-p)/{p}}\,\left(\sqrt{-\frac{1}{2\sigma^2\pi\alpha}}\right)^n\,\left(\sqrt{-\frac{\alpha}{2\sigma^2}}\right)^
{n\frac{1-p}{2p}}\, \Gamma\Big(n(p-1)/{2p}\Big)\,+O(\varepsilon)}=$$
$$
=\,-\,C\,\varepsilon\,x\,|x|^{\frac{2(1-p)}{p}}+o(\varepsilon),\quad
\varepsilon\to 0,\eqno(4.3)$$ with the positive constant
$$C=\,-\frac{1}{\alpha}\,\left(\sqrt{-\frac{\alpha}{2\sigma^2}}\right)^{\frac{1-p}{2p}}\,
\frac{\Gamma\Big((n+2)(p-1)/2p\Big)}{\Gamma\Big(n(p-1)/2p\Big)}.\eqno(4.4)
$$

 If $p\le 1,$ the integrals
in (4.1) and (4.2) diverge as $L\to \infty.$ Now we take into
account that
$$
\int\limits_{\Omega_L}\,P(t(\varepsilon),x,u)\,du=\,f_L\,\int\limits_0^{|\alpha
L|}\,\sum\limits_{k=0}^{\infty}\,F_k(x,|u|)\,\varepsilon^k\,d|u|,$$
$$F_k(x,|u|)=e^{\frac{\alpha\,|x|^2}{2\sigma^2
|u|^{2p}}}\,\sum\limits_{s=0}^{k}\frac{c_{ks}\,|x|^{2s}}{|u|^{sp+np-n+1}},$$
$$
\int\limits_{\Omega_L} \,u P(t(\varepsilon),x,u)\,d u=\varepsilon
x\,\int\limits_0^{|\alpha
L|}\,\sum\limits_{k=0}^{\infty}\,G_k(x,|u|)\,\varepsilon^k\,d|u|,$$
$$
G_k(x,|u|)=e^{-\frac{\alpha\,|x|^2}{2\sigma^2
|u|^{2p}}}\,\sum\limits_{s=0}^{k}\frac{b_{ks}\,|x|^{2s}}{|u|^{(s+1)p+np-n-1}},$$
with some constant $ b_{ks}\,$ and $\, c_{ks}.$ Further,
$$\hat u (t(\varepsilon)),x)\,=\,\varepsilon
x\,\lim\limits_{L\to\infty}\,\frac{\int \limits_0^{|\alpha
L|}\,\sum\limits_{k=0}^{\infty}\,G_k (x,\xi)\,\varepsilon^k\,d\xi}
{\int \limits_0^{|\alpha
L|}\,\sum\limits_{k=0}^{\infty}\,F_k(x,\xi)\,\varepsilon^k\,
\,d\xi}=
$$
$$=\,x\,\varepsilon\,\lim\limits_{L\to\infty}\,\frac
{\frac{\omega_n}{\alpha}\left(\sqrt{-\frac{1}
{2\sigma^2\pi\alpha}}\right)^n \,\int \limits_0^{|\alpha L|}\,
e^{\frac{\alpha |x|^2}{2\sigma^2
\xi^{2p}}}\,\frac{1}{\xi^{(n+2)p-n-1}}\,d\xi\,(1+O(\varepsilon))
}{\omega_n\,\left(\sqrt{-\frac{1}{2\sigma^2\pi\alpha}}\right)^n
\,\int \limits_0^{|\alpha L|}\, e^{\frac{\alpha |x|^2}{2\sigma^2
\xi^{2p}}}\,\frac{1}{\xi^{n
p-n+1}}\,d\xi\,(1+O(\varepsilon))}=\eqno(4.5)
$$
$$
=\,\left\{\begin{array}{ll}
\frac{x\,\varepsilon}{\alpha}+o(\varepsilon),&\quad p=1, \\ -\,
\infty\,\sign x,&\quad p<1.
\end{array}\right.,
$$
where we use de L'H$\rm \hat{o}$pital's rule to compute the ratio of
the divergent integrals in (4.5).

Thus, for $p=1$ we  get (4.3) again.

For $p<1$ we are going to
obtain a more specified result. 
Recall that
$$ P(t(\varepsilon),x,u)\,
=\eqno(4.6)$$
$$=\,f\left(\frac{u}{\alpha}\right)\,\left(
\,\frac{1}{\left(\frac{2\pi\sigma^2(\varepsilon-1)}{\alpha}\right)^{n/2}}\,
\frac{\exp\left(-\frac{|u|^{2-2p}\,\varepsilon^2}{2\alpha\sigma^2
(\varepsilon-1)}\right)}{|u|^{np}}\right)\,\left( \exp\left(
-\frac{\alpha(|x|^2-2(x,u)\varepsilon/\alpha)}{2\sigma^2
(\varepsilon-1)|u|^{2p}}\right)\right).$$ For  the third factor as
$\varepsilon\to 0$ we have:
$$e^{{\frac
{2\alpha {|x|}^{2}}{{\sigma^2}\,{|u|}^{2p}}}}+e^{{\frac {2\alpha{
|x|}^{2}}{{\sigma^2}\,{|u|}^{2p}}}} \left( -\frac
{(x,u)}{{\sigma^2}\,{|u|}^{2p}} \, + \,\frac {\alpha
{|x|}^{2}}{{2\sigma^2}\,{|u|}^{2p}} \right) \,{\varepsilon}+O
 \left( {{\varepsilon}}^{2} \right) .$$
The second factor in (4.6) secures the convergence of integrals.
Calculations of both integrals in (2.2) with the use of the Maple
environment allow to obtain explicit formulas for rational $p,$
the result being expressed through special functions (Bessel, Gamma
and hypergeometric functions). We do not quote here this formula as
it is  very cumbersome.  The simplest result is for $p=\frac{1}{2}:$
$$\lim_{L\to\infty}\,\int\limits_{\Omega_L}\,P(x,u,t)\,du \,=
2 \,f_L\,K\left(\frac{n}{2},\frac{\varepsilon |x|}{\sigma^2
\sqrt{1-\varepsilon}}\right)\,\left(\frac{|x|}{2\sigma^2\pi\varepsilon\sqrt{1-\varepsilon}}\right)^{\frac{n}{2}}\,(1+O(\varepsilon))$$
$$\lim_{L\to\infty}\,\int\limits_{\Omega_L}\,u\,P(x,u,t)\,du
\,=$$
$$
 =\,f_L\,\frac{2  \alpha \sqrt{1-\varepsilon} x}{ n \sigma^2}\,K\left(\frac{n}{2}+1,\frac{\varepsilon
 |x|}{\sigma^2
\sqrt{1-\varepsilon}}\right)\,\left(\frac{|x|}{2\sigma^2\pi\varepsilon
\sqrt{1-\varepsilon}}\right)^{\frac{n}{2}}\,(1+O(\varepsilon)),
$$
where $K$ is the modified Bessel function of the second kind
\cite{Ryzhik}. Thus,
$$\hat u(t,x)\,=\,\frac{2\alpha  \sqrt{1-\varepsilon} x |x|}{n \sigma^2}\,\frac{K\left(\frac{n}{2}+1,
\frac{\varepsilon |x|}{\sigma^2
\sqrt{1-\varepsilon}}\right)}{K\left(\frac{n}{2},\frac{\varepsilon
|x|}{\sigma^2 \sqrt{1-\varepsilon}}\right) }\,(1+O(\varepsilon)).
$$
Calculations using asymptotic expansion of the $K$ functions for
$\varepsilon \to 0$ show that
$$\hat u(t,x)\,= \,\frac{\alpha}{\varepsilon}\,x\,+\,O\left(1\right),\quad \varepsilon\to 0.$$

Analogously result we get in the case of any rational $p<1.$ Namely,
if we assume $p=\frac{m_1}{m_2}, \quad m_1, m_2 \in {\mathbb N}, \,$
we obtain
$$\hat u (t,x)\,= \,\frac{\alpha}{\varepsilon}\,x\,+\,o\left(\frac{1}{\varepsilon} \right),\quad \varepsilon\to 0.
\eqno(4.7)$$ In particular, for $p=\frac{1}{m},\quad m\in {\mathbb
N}$
$$\hat u(t,x) \,=\, \frac{\alpha}{\varepsilon}\,x\,+ O\left(\varepsilon^{-\frac{2m-3}{2m-1}} \right),\quad \varepsilon\to 0.$$

Then, we can consider $p$ as a parameter and notice that at any
fixed $x\ne 0, \,\varepsilon \in (0,1]$ the functions
$P(t(\varepsilon), x,u;p)$ and $u P(t(\varepsilon), x,u;p)$ are
continuous on the set $u\in\Omega_L,\,p\in [0,1].$ Moreover,
$\frac{1}{f_L}\,\lim\limits_{L\to\infty}\,
\int\limits_{\Omega_L}\,P(t(\varepsilon), x,u;p)\,du $ and
$\frac{1}{f_L}\,\lim\limits_{L\to\infty}\,
\int\limits_{\Omega_L}\,u\,P(t(\varepsilon), x,u;p)\,du\,\,$ are
uniformly bounded for $p\in [0,1-\delta],\,$ $ \delta$ is a positive
arbitrary small constant (see (4.6)). Thus, by a standard reasoning
we conclude that the ratio $\hat u(t,x)$ is continuous in $p,\,p\in
[0,1-\delta].$ Since there exist a rational sequence convergent to
every real $p,$ the continuity implies that (4.7) holds for all
$p\in [0,1-\delta].$ Over the arbitrariness of $\delta$ we can
conclude that (4.7) is true for all $p\in [0,1).$

\medskip

Now we fix $\varepsilon>0$ and expand the integrals in (2.2) near
$x=0.$ We have for $\beta=0,$ any $k=1,...,n,$ using (2.9):
$$
\int\limits_{\Omega_L}\,u_k\,P(t(\varepsilon),x,u)\,du \,=$$$$=\,
\frac{f_L\,\varepsilon\,x_k}{ \sigma^2
(\varepsilon-1)\left(2\pi\sigma^2(\varepsilon-1)\alpha
\right)^{n/2}}\,\omega_n\,\times$$$$\times \int\limits_0^{|\alpha
L|}\,\exp\left(-\frac{|u|^{2-2p}\,\varepsilon^2}{2\alpha\sigma^2
(\varepsilon-1)}\right)\,|u|^{n(1-p)-1-2p}\,u_k^2\,d|u|+o(|x|^2)=$$$$=
c_k(\varepsilon)\,x_k,
$$
for some functions $c_k(\varepsilon).$ Moreover,
$$
\int\limits_{\Omega_L}\,P(t(\varepsilon),x,u)\,du \,=$$$$\,
\frac{f_L\,\omega_n}{\left({2\pi\sigma^2(\varepsilon-1)\alpha}\right)^{n/2}}\,
\int\limits_0^{|\alpha
L|}\,\exp\left(-\frac{|u|^{2-2p}\,\varepsilon^2}{2\alpha\sigma^2
(\varepsilon-1)}\right)\,|u|^{n(1-p)-1}\,d|u|+o(|x|)=
$$
$$
=\frac{f_L}{2|p-1|}\,\left(\varepsilon^2\,\pi\right)^{-\frac{n}{2}}\,\Gamma\left(\frac{n}{2}\right)+o(|x|).
\eqno(4.8)
$$
All integrals converge for $\,p>0,\,p\ne 1$ as $L\to\infty.$
Further,
$$A(\varepsilon):=\sum\limits_{j=k}^n\,c_k(\varepsilon)=\frac{\alpha\,}{\varepsilon |p-1|}\,\left(\varepsilon^2\,\pi\right)^{-\frac{n}{2}}\,
\Gamma\left(1+\frac{n}{2}\right),$$ due to the equality of space
directions
$c_k(\varepsilon)=\frac{1}{n}\,A(\varepsilon),
$ therefore
$$\hat u(t,x)=\,\frac{\alpha
}{\varepsilon}\, x\,+\,o(|x|),\quad |x|\to 0.\eqno(4.9)$$ In the
limit case $p=1,$ where both integrals in the ratio  (2.2) diverge,
we can apply de L'H$\rm \hat{o}$pital's rule as above (in (4.9) to
obtain (4.6)).

Thus, the theorem is proved.

\end{proof}

\begin{remark}
For $\beta>0,\,p\ge 1 $ we get the expansion (4.3) as
$\varepsilon\to 0$ again, but in this case
$$
C= -\frac{1}{\beta}\,\ln \frac
{\alpha+\beta}{\alpha}\,\left(2\,\left( \left( {\frac
{\alpha}{\beta+\alpha}} \right) ^{2p} -1 \right) \frac{\sigma}{2p
\beta}\right)^{\frac{p-1}{p}}\,\frac{\Gamma\Big((n+2)(p-1)/2p\Big)}{\Gamma\Big(n(p-1)/2p\Big)},
$$
for $p\in [0,1)$ we get qualitatively the same result as in (4.7).

 Analogously we have for any $p>0, \,\varepsilon\in
(0,1]$
$$\hat u(t,x)=\,\frac{\beta}{\left(\frac{\beta}{\alpha}+1 \right)^\varepsilon-1}\,\,x\,+\,o(|x|),\quad |x|\to 0.
$$

\end{remark}

\medskip
\medskip

\begin{remark}
 Let us compute the conditional variance $\hat v (t,x).$ 
For $p>1+\frac{4}{n}$ we can expand the probability density
$P(t,u,x)$ in $\varepsilon$ near $\varepsilon=0$ and obtain that for
any $x\ne 0$

$$\hat v (t(\varepsilon),x)=
F(n,p,\alpha,\sigma)\,|x|^{\frac{2}{p}}\,\left(1-\frac{n(p-1)-(p+2)}{2p}\,\varepsilon\right)+o(\varepsilon),
$$
with
$$
F(n,p,\alpha,\sigma)= \frac{\Gamma\left(
\frac{n(p-1)-2}{2p}\right)}{\Gamma\left(
\frac{n(p-1)}{2p}\right)}\,\left(-\frac{\alpha}{4\sigma^2}\right)^{\frac{1}{p}}.
$$
Thus, at any fixed $x\ne 0$ the variance tends to some finite value
as $\varepsilon\to 0.$

Further, for any fixed $\varepsilon>0$ we get the following
asymptotic expansion near $x=0:$
$$\hat v (t(\varepsilon),x)
=\frac{\Gamma\left( \frac{n(p-1)-2}{2p}\right)}{\Gamma\left(
\frac{1}{2}n\right)}\,\,\left(-\frac{\varepsilon^2}{4\,\alpha\,\sigma^2\,(1-\varepsilon)}\right)^{\frac{1}{p-1}}+o(|x|),
$$
therefore the variance tends to zero as $x\to 0$ and $\varepsilon\to
0.$

We cannot write an explicit formula for the variance for all $p\le
1+\frac{4}{n}.$ However, for $p\le 1$ in the expression for the
second moment both integrals in the ratio diverge as $L \to \infty,$
and we can use de L'H$\rm \hat{o}$pital's rule once more to show
that the  variance tends to infinity as $\varepsilon \to 0.$ In
particular, for $p=\frac{1}{m},\quad m\in{\mathbb N},$
$$\hat v (t,x)
=O\left(\varepsilon^{-\frac{4m}{2m-1}}\right),\quad \varepsilon \to
0.$$

\end{remark}

\section{Asymptotic behavior for a Gaussian initial distribution}

Let us set $f=\left(\frac{k}{\sqrt{\pi}}\right)^n\,\exp(-k^2
x^2),\quad k>0,$ and study the asymptotics of (2.2), (2.9) near
$x=0.$ We consider only the case $\beta=0$ and $n=1.$ As we have
seen in Sec.3, this type of initial density distribution eliminates
the unbounded gradient growth as $\varepsilon\to 0$ in the case of a
constant $\sigma.$

To get a qualitative result we firstly find an asymptotic expansion
near $x=0$ for both integrals in(2.2), and for their ratio in
$\varepsilon$ at the point $\varepsilon =0\quad (t=T).$ The
computations made in the environment Maple show that for $p<1$
$$
\hat u=-C\,\varepsilon \,x\,+\, o(\varepsilon) C_1(x)\,+\,o(|x|)
C_2(\varepsilon),
$$
and for $p\ge 1$
$$
\hat u=-\frac{\alpha}{\varepsilon}\,x \,+\,C\,\,x\,+\,
o(\varepsilon) C_1(x)\,+\,o(|x|) C_2(\varepsilon),
$$
with a positive constant $C$ depending only on $\,\alpha, \,\sigma,
k\,$ and bounded functions $C_1$ and $C_2.$

Numerical computations show that if $p < 1$ the function $\hat u
(t,x)$ is very close to the respective function for $p=0.$ In
particular, the singularity at the origin is also eliminated. If
$p\ge 1,$ the picture is similar to the case of uniform initial
distribution $f(x).$

\section{"Observable" density and "induced" velocity}

Introduce the function
$$\hat \rho(t,x)=\lim_{L\to\infty}\,\int\limits_{\Omega_L}\,P(t,x,u)\,du.$$
If $\sigma=0,$ in (2.1), then $\hat \rho (t,x)$ solves the
continuity equation
$$
\frac{\partial \hat \rho}{\partial t}\,+\,{\rm div} (u\,\hat
\rho)=0,\eqno(6.1)
$$
and therefore corresponds to the distribution of density for the
particles provided initial data $\hat \rho(0,x)=f(x)> 0$ are given.
Here we do not require that $f(x)$ is a probability density. Thus,
$\hat\rho(t,x) $ coincides with the real density of particles that
we denote by $\rho(t,x).$

For $\sigma\ne 0$ the function $\hat \rho(t,x)$ is not a solution to
(6.1). Nevertheless an observer who does not know the random
distribution of the particles positions considers $\hat \rho (t,x)$
as the density for the particles.  We will call $\hat \rho(t,x)$ in
this case an observable density. It differs from $\rho (t,x)$. For
example, it is easy to compute that for $f(x)=C$  and $\beta=0$ from
the system (1.1), (6.1) one gets $\rho(t,x)\,=\,\frac{C}{(1+\alpha
t)^n}\,=\,C\,\varepsilon^{-n}.$ The observable density $\hat
\rho(t,x),$ as follows from (4.8), behaves in a quite different way.
The only common feature is the asymptotics $O(\varepsilon^{-n})$ as
$|x|\to 0.$ Nevertheless, given a density $\rho$, one can find from
(6.1) the velocity $u(t,x).$ Let us consider the case $n=1.$ If we
assume
$$\lim\limits_{|x|\to \infty} \rho(t,x) u(t,x) =0,\quad t\in \mathbb R_+ $$ (this
implies the momentum conservation), we get
$$
u(t,x)=-\frac{\int\limits_{-\infty}^{x}\,\rho'_t(t,x)\,dx}{\rho(t,x)}.
$$
Remark that $\hat\rho(t,x)$ is positive for $f(x)>0.$ Let us
introduce the vector-function $v(t,x)$ according to the formula
$$
v(t,x)=-\frac{\int\limits_{-\infty}^{x}\,\hat
\rho'_t(t,x)\,dx}{\hat\rho(t,x)}.\eqno(6.2)
$$
It is natural to call $v(t,x)$ the "induced" velocity. As follows
from (6.2), (2.1), we have
$$v(t,x)= \hat u(t,x)+ v_1(t,x),\eqno(6.3)$$
where $$v_1(t,x):= \frac {\sigma^2}{2}\,\frac{\int\limits_{{\mathbb
R}}\,|u|^{2p}\,P'_x(t,x,u)\,du} {\int\limits_{{\mathbb
R}}\,P(t,x,u)\,du}.$$ One can see that if $P(t,x,u)\,$  is given by
(2.9) in the case $f(x)=const$ the second term $v_1$ in (6.3)
vanishes at $x=0.$ Moreover, computations show that as $|x|\to 0$
$$v_1(t,x)\,=\,C(t)\,\varepsilon \,x \,|x|^2\,+\,o(|x|^3),$$
with  a function $C(t)$  which stays bounded as $\varepsilon \to 0.
$ Thus, in the vicinity of the point of the singularity formation
the induced velocity is close to  $\hat u.$

\section{Conclusion}

We have shown that in the model (1.8), generally speaking, the
expectation of the velocity given position differs drastically from
the velocity in the Burgers equation in the deterministic coordinate
space. Under the observation in the random position of particles the
blow up phenomena can be lost. In this context the uniform initial
distribution of the particles position seems most interesting.  Here
the threshold effect arises. Namely, if the exponent $p$ is less or
equal then 1, then the expectation $\hat u$ follows the real
behavior of the velocity rather well. However, if $p\ge 1,$ the
function $\hat u$ tends to zero whereas the real velocity tends to
infinity at any point outside the origin $x=0$ as the time tends to
the critical value.

The question on a PDE obeyed by the function $\hat u$ arises
naturally, due to existing formalisms to represent solutions of PDE
as the expected value of functionals of stochastic processes (see
e.g. \cite{Freidlin},\cite{Freidlin1}, \cite{Belopolskaya} and
references therein). In \cite{Iyer} one can find a recent result in
the field close to our study, the stochastic formulation of the
viscous Burgers equation. Namely, it was shown that if the pair
$(u,X)$ (where $X$ is the flow map such that $u\circ X$ is constant
in time) solve the stochastic system
$$
dX=\,u\,dt + \sqrt{2\nu}\,dW,\quad u=\textbf{E}[u_0\circ
X^{-1}],\quad X(x_0,0)=x_0,\quad \nu=const>0,
$$
then under certain additional conditions $u$ satisfies the viscous
Burgers equation
$$
u_t+(u,\nabla)\,u=\,\nu\,\Delta u.\eqno(7.1)
$$
The function $\hat u$ also solve (7.1) for $p=0,$ however only for
the special initial data and uniform initial distribution of
particles. For example, one can check that $\hat u$ defined as
(2.2), with $P(t,x,u)$ governed by (2.1) does not satisfies the
viscous Burgers equation even for $p=0$ and linear initial data for
the Gaussian initial distribution (see (3.5)). Nevertheless, the
comparison of $\hat u$ and the solution to the Burgers equation with
a specific viscosity term is a very interesting open question.

At last we would like to discuss possible applications of the
results of this paper. First of all the the stochastic systems like
(1.8) with velocity field of the fluid, with possibly random
component with prescribed statistics can be applied to model
turbulent or other disordered fluctuations. For example, the problem
can be in describing some desired statistics of the trajectory
$X(t)$ of a tracer particle released initially from some point $x_0$
and subsequently transported jointly by the flow $u(x, t)$ and
molecular diffusivity \cite{Majda}.

Other applications can be found in the theoretical financial
mathematics to model the earnings yield of a risky asset ($X(t)$)
with the trend and volatility depending on a certain macroeconomic
factor $U(t)$ (e.g. the spot interest rate) in the spirit of
\cite{BP}. Let us notice in this context that for $p=0$ in several
cases it is possible to solve the Fokker-Planck equation for more
general class of stochastic DE, describing the economic quantities
in a more agreed way, namely
$$
d X_k(t)\,=\,(a_k+A_{ki} U_i(t))\,dt\,
+\,\sum\limits_{i=1}^{m+n}\sigma_{ki} \,d(W_i)_t,\quad k=1,..,n;
$$
$$d U_i(t)=\,(b_i+B_{ik}\,U_k(t))\,dt\,+\,\sum\limits_{k=1}^{m+n}\lambda_{ik} \,d(W_k)_t   ,\quad i=1,..,m;
$$
$$
X_k(0)=x_k,\quad U_i(0)=u_i,\quad t\ge 0,
$$
where $W(t)$ is a ${\mathbb R}^{m+n}$ valued standard Brownian
motion process, the market parameters $a_k, A_{ik}, \sigma_{ik},
b_k,\,B_{ik}, \lambda_{ik}$ are constant matrices of appropriate
dimensions.

\section*{Acknowlegments}

Stimulating discussion with Y.I.Belopolskaya and V.M.Shelkovich are
gratefully acknowledged. We thank M.Freidlin and L.Ryzhik for
attracting attention to the paper \cite{Iyer}.



\begin{thebibliography}{99}
\bibitem{Hopf} E.Hopf, { The partial differential equation $u_t+uu_x=\mu u_{xx},$}
{\it Comm.Pure Appl.Math} {\bf 3}\,(1950)  201-230.
\bibitem{AKl1} S.Albeverio, A.Klar,{
Long time behavior of nonlinear stochastic oscillators: The one-
dimensional Hamiltonian case.} {\it J. Math. Phys.} {\bf 35} (8)\,
(1994) 4005-4027.

\bibitem{AK2} S.Albeverio, A. Klar,
{Longtime behaviour of stochastic Hamiltonian systems: the
multidimensional case.} {\it Potential Anal.} {\bf 12} (3) (2000)
281-297.
\bibitem{AK} S.Albeverio, V.N. Kolokoltsov,
{ The rate of escape for some Gaussian processes and the scattering
theory for their small perturbations.} {\it Stochastic Processes and
their Applications} {\bf 67} (2) (1997) 139-159.
\bibitem{AHZ}
S. Albeverio, A. Hilbert, E. Zehnder, { Hamiltonian systems with a
stochastic force: nonlinear versus linear and a Girsanov formula.}
{\it Stochastics and Stochastics Reports} {\bf 39} \,(1992),
159-188.

\bibitem{Risken}H. Risken, {\it The Fokker-Planck
Equation Methods of Solution and Applications,} Second Edition,
(Springer-Verlag, 1989).
\bibitem{Woy} W.A. Woyczy$\rm\acute{n}$ski, {\it Burgers-KPZ Turbulence},
LNM 1700, Springer, 1998.
\bibitem{Khanin} J. Bec, K.Khanin, { Burgers turbulence,}
{\it Physics Reports}  {\bf 447} (1)\,(2007), 1 - 66 .
\bibitem{Bouchaud}J.P.Bouchaud, M.Mezard, { Velocity fluctuations in forced
Burgers turbulence,}{\it Phys.Rev. E} {\bf 54} (1996)\, 5116
\bibitem{Gurarie}V.Gurarie, { Burgers equations revisited,}
{\it arXiv:nlin/0307033v1\, [nln.CD]}\, (2003).
\bibitem{oksendal} B.{\O}ksendal, {\it Stochastic differential equations: an introduction
with applications,} 5th ed.(Springer-Verlag Heidelberg New York,
1998).
\bibitem{Chorin} A.J.Chorin, O.H.Hald, {\it Stochastic tools
in mathematics and science}  (New York: Springer, 2006).

\bibitem{Ryzhik} I.S. Gradshteyn, I.M. Ryzhik,
{\it Table of integrals, series, and products,}  6th ed.  (San
Diego, CA:
Academic Press, 2000).
\bibitem{Freidlin}M.Freidlin,
{\it  Functional integration and partial differential equations},
Annals of Mathematics Studies, No.109 (Princeton, New Jersey:
Princeton University Press, 1985).

\bibitem {Freidlin1} M.Freidlin, {\it Markov processes and differential equations:
asymptotic problems}, Lectures in Mathematics, (ETH Zu"rich. Basel:
Birkha"user, 1996).

\bibitem{Belopolskaya} Ya.I.Belopolskaya; Yu.L.Daletskij,
{\it Stochastic equations and differential geometry} (Kluwer
Academic Publishers, 1990).

\bibitem{Iyer} P.Constantin, G.Iyer,\, A stochastic Lagrangian
representation of the three-dimensional incompressible Navier-Stokes
equations. {\it Commun. Pure Appl. Math.} {\bf 61} (3) 330-345
(2008).


\bibitem{Majda} A.J. Majda, P.R. Kramer,{ Simplified models for turbulent diffusion: Theory, numerical
modelling, and physical phenomena} {\it Physics Reports} {\bf 314}
237-574 (1999).

\bibitem{BP} T.R.Bielecki, S.R.Pliska,  Risk-sensitive
dynamic asset management. {\it Appl. Math. Optimization} {\bf 39}
(3) 337-360 (1999).
\end{thebibliography}
\end{document}